\documentclass[12pt]{article}

\usepackage[latin1]{inputenc}
\usepackage[T1]{fontenc}
\usepackage{kpfonts}
\usepackage[french,english]{babel}
\usepackage{xspace}
\usepackage{eurosym} %symbole euro
\usepackage{enumerate} %pour modifier les compteurs des enumerations
\usepackage{verbatim} %mise en page ignoree
\usepackage{cancel}
\usepackage{changepage}
\usepackage[pdftex,leftbars,xcolor]{changebar}
\usepackage{multicol} %texte sur plusieurs colonnes

\usepackage{multirow}
\usepackage{tabularx}

\allowdisplaybreaks[2]

\usepackage[pdftex,a4paper,centering]{geometry}
\usepackage{pdflscape}

\usepackage{fancyhdr}

\usepackage{amsmath,amsthm}
\usepackage{amssymb,latexsym,amsfonts}

\usepackage{textcomp}

\usepackage{aliascnt}
\numberwithin{equation}{section}

\usepackage{mathtools}
\usepackage{relsize,exscale}
\usepackage{rotating}

\usepackage{etoolbox}

\usepackage{etex}

\usepackage[all]{xy}

\usepackage{xcolor}

\pdfoutput=1
\usepackage[pdftex,%
bookmarks=true,%
bookmarksnumbered=true,%
plainpages=false,%
breaklinks=true,%
colorlinks=false,%
urlcolor=blue,%
citecolor=green,%
linkcolor=red,%
pdftitle={Heisenberg formality},%
pdfauthor={Olivier Elchinger}]{hyperref}
\usepackage{doi}

\usepackage[hyperpageref]{backref}
\backreffrench
\renewcommand*{\backref}[1]{}  % Disable standard
\renewcommand*{\backrefalt}[4]{% Detailed backref
  \ifcase #1 %
  \relax%(Not cited.)%
  \or
(Cited page~#2.)%
  \else
(Cited pages~#2.)%
  \fi}

\newaliascnt{cor}{theorem}

\aliascntresetthe{cor}

\newaliascnt{prop}{theorem}

\aliascntresetthe{prop}

\newaliascnt{lemma}{theorem}

\aliascntresetthe{lemma}

\theoremstyle{definition}
\newaliascnt{defi}{theorem}
\newtheorem{defi}[defi]{Definition}
\aliascntresetthe{defi}

\newaliascnt{example}{theorem}
\newtheorem{example}[example]{Example}
\aliascntresetthe{example}

\theoremstyle{remark}
\newaliascnt{remark}{theorem}

\aliascntresetthe{remark}

\addto\extrasenglish{%
}

\newcommand{\ie}{\textit{i.e.} \/}

\newcommand{\numberset}[1]{\mathbb{#1}}

\newcommand{\nat}{\numberset{N}}

\newcommand{\korps}{\numberset{K}}

\DeclareMathOperator{\Kr}{Ker}
\DeclareMathOperator{\Img}{Im}

\DeclareMathOperator{\Hom}{Hom}

\DeclareMathOperator{\Sym}{\mathcal{S}\!}

\newcommand{\ChEcochains}[1][]{\ensuremath{C_{CE}^{#1}}}
\newcommand{\ChE}[1][]{\ensuremath{H_{CE}^{#1}}}

\title{Study of formality for the Heisenberg algebra}
\author{Olivier Elchinger \\
Laboratoire de Math\'{e}matiques, Informatique  et Applications \\
Universit\'{e} de Haute Alsace, Mulhouse, France
}

\begin{document}

\maketitle

%\vspace*{2em}

\begin{abstract}
In this paper, we compute the Chevalley-Eilenberg cohomology of the three-dimensional Heisenberg Lie algebra with values in its universal enveloping algebra. We also compute the Schouten brackets on cochains and cohomology level in order to write the formality equations. It turns out that there is no formality, and that the perturbed $L_\infty$ structure on the cohomology has non-trivial terms in infinitely many degrees.
\end{abstract}



\section*{Introduction}

Deformation of structure theories is useful to formalize quantum physics. The algebraic structure considered by classical mechanics is the associative commutative algebra of smooth functions over a symplectic, or more generally, a Poisson manifold. Deformation quantization consists to construct an associative non-commutative multiplication (more precisely a $\star$-product) on the formal series in $\hbar$ (which represents Planck's constant) with coefficients in this algebra, which encodes the Poisson bracket in the first order. The Poisson structure is then called the classical limit and the deformation is the $\star$-product. This point of view, initiated in 1978 in \cite{BFFLS78}, tries to consider quantum mechanics as a deformation from classical mechanics. This pioneering work raises the fundamental questions about the existence and the uniqueness of a deformation quantization for a given Poisson manifold.

The first results treated the case of symplectics manifolds. The general case of Poisson manifolds was solved by Kontsevich in 1997 in \cite{K03}. He deduced the result by proving a much more general statement, which he called \og{}formality conjecture\fg{}. Endowed with the Gerstenhaber bracket, the Hochschild complex of the algebra of smooth functions over a Poisson manifold admits a graded Lie algebraic structure by shift, which controls the deformations of the Poisson bracket. Kontsevich shows that this complex is linked with its cohomology --- which therefore controls the same deformations --- by a $L_\infty$-quasi-isomorphism, called a formality map. Using this formality map, he shows of to deform the commutative product of functions in a $\star$-product.

Konstevich shows in particular that the formality criterion is true for symmetric algebras over a finite-dimensional vector space are formal. Bordemann and Makhlouf have examined in \cite{BM08} a slight generalization to the case of universal enveloping algebras. They showed (implicit in Kontsevich's work) that the formality for a Lie algebra is equivalent to those of its universal enveloping algebra. They also proved that there is formality for the universal enveloping algebra of an affine Lie algebra. These methods were used in \cite{BMP05} and give informations about the rigidity of universal enveloping algebras. \\

In \cite{El12}, we detailed the definitions of formality for associative and Lie algebras, and we studied the formality equations for the associative free algebra, and for the Lie algebra $\mathfrak{so}(3)$. It turns out that in these cases, there is no formality in Kontsevich's sense, but one has to add an extra term to the graded Lie bracket on the cohomology to obtain a perturbed formality map.

This work is part of this line of research. We study here the formality for the three-dimensional Heisenberg Lie algebra. Although the Lie bracket is more simple than the one of $\mathfrak{so}(3)$, the cohomology is more important, and therefore the formality equations contains more terms. Actually, to obtain a perturbed formality map in this case, we have to add infinitely many terms to the Schouten bracket on the cohomology. \\

The paper is organised as follows. In \autoref{Sec:formality-Lie}, we set the notations and recall the graded Lie structure (Schouten bracket) on the Chevalley-Eilenberg complex of a Lie algebra and its cohomology. We also give the definition of formality for this complex and write down the formality equations. In \autoref{Sec:Heisenberg-cohomology}, we compute the cohomology of the Heisenberg algebra. \autoref{Schouten-cochain} and \autoref{Schouten-cohomology} are devoted to the computation of the Schouten brackets on cochain level and cohomology, in order to solve the formality equations in \autoref{Sec:formality-eq}. \\

I wish to thank gratefully Mr Makhlouf who showed me his first computations of the cohomology for the Heisenberg algebra and for his advice to pursue the study of the formality.

\section{Formality of Lie algebras} \label{Sec:formality-Lie}

Let $(\mathfrak{g},[~,~])$ be a Lie $\korps$-algebra. We consider the Chevalley-Eilenberg differential complex $(\ChEcochains(\mathfrak{g},\Sym\mathfrak{g}),\delta_{CE})$ with values in $\Sym\mathfrak{g}$
\begin{gather}
C^0 \coloneqq \ChEcochains[0](\mathfrak{g},\Sym\mathfrak{g}) \coloneqq \Sym\mathfrak{g} \qquad C^k \coloneqq \ChEcochains[k](\mathfrak{g},\Sym\mathfrak{g}) \coloneqq \Hom\big({\bigwedge}^k \mathfrak{g},\Sym\mathfrak{g}\big) \quad \text{for $k \in \nat^\star$} \\
\delta_{CE}^k : C^k \to C^{k+1} \notag
\end{gather}
and its cohomology
\( \displaystyle
\ChE(\mathfrak{g},\Sym\mathfrak{g}) \coloneqq \bigoplus_{k \in \nat} \ChE[k](\mathfrak{g},\Sym\mathfrak{g})
\)
where for $k \in \nat$, the $k^\mathrm{th}$ Chevalley-Eilenberg cohomology group with coefficients in $\Sym\mathfrak{g}$ is
\begin{equation}
H^k \coloneqq \ChE[k](\mathfrak{g},\Sym\mathfrak{g}) \coloneqq \frac{\Kr(\delta_{CE}^k)}{\Img(\delta_{CE}^{k-1})} \coloneqq Z^k/B^k
\end{equation}

The coboundary operator $\delta_{CE}$ equals $[\pi,~]_s$, where $\pi = [~,~] \in \Hom\big({\bigwedge}^2 \mathfrak{g},\mathfrak{g}\big)$ is the Lie bracket from $\mathfrak{g}$ and $[~,~]_s$ is the Schouten bracket on the complex $\ChEcochains(\mathfrak{g},\Sym\mathfrak{g})$.

As the Gerstenhaber bracket, the Schouten bracket defines a graded Lie bracket on the shifted space $\ChEcochains(\mathfrak{g},\Sym\mathfrak{g})[1]$. Therefore the graded antisymmetry and graded Jacobi identities: for $\xi,\eta,\zeta \in \ChEcochains(\mathfrak{g},\Sym\mathfrak{g})[1]$ we have
\begin{gather}
[\eta,\xi]_s = -(-1)^{|\eta| |\xi|}[\xi,\eta]_s, \label{gr-antisymmetry} \\
(-1)^{|\xi| |\zeta|} [[\xi,\eta]_s,\zeta]_s + (-1)^{|\eta| |\xi|} [[\eta,\zeta]_s,\xi]_s + (-1)^{|\zeta| |\eta|} [[\zeta,\xi]_s,\eta]_s = 0 \label{gr-Jacobi},
\end{gather}
degrees being taken in $\ChEcochains(\mathfrak{g},\Sym\mathfrak{g})[1]$.

On the other hand there is the pointwise exterior multiplication $\wedge$ which makes $\ChEcochains(\mathfrak{g},\Sym\mathfrak{g})$ a graded commutative algebra. The Schouten bracket and the exterior multiplication are compatible by the graded Leibniz rule
\begin{equation}
[\xi,\eta \wedge \zeta]_s = [\xi,\eta]_s \wedge \zeta + (-1)^{|\xi| (|\eta|+1)} \eta \wedge  [\xi,\zeta]_s, \label{gr-Leibniz}
\end{equation}
(if $\eta$ is of degree $|\eta|$ in $\ChEcochains(\mathfrak{g},\Sym\mathfrak{g})[1]$ then it is of degree $|\eta| + 1$ in $\ChEcochains(\mathfrak{g},\Sym\mathfrak{g})$). \\

Throughout the computations, we will use the fact that $\ChEcochains(\mathfrak{g},\Sym\mathfrak{g})$ is an $\Sym \mathfrak{g}$-module \ie for $\xi \in \ChEcochains(\mathfrak{g},\Sym\mathfrak{g}),\ f \in \Sym\mathfrak{g}$, $f \wedge \xi = f \xi$. We also give some examples of Schouten brackets which will be used later.
\begin{example}
The space $C^1 \coloneqq \ChEcochains[1](\mathfrak{g},\Sym\mathfrak{g})$ is the space of vectors fields, \ie the derivations of $\Sym\mathfrak{g} = C^0$. For $f,g \in \Sym\mathfrak{g},\ X,Y \in C^1$
\begin{gather*}
[f,g]_s = 0, \\
[X,f]_s = X(f) = \sum_{i=1}^n X_i \partial_i f, \\
[X,Y]_s = [X,Y] = \sum_{i=1}^n X_i \partial_i Y - \sum_{i=1}^n Y_i \partial_i X,
\shortintertext{with in particular}
[f \partial_i,g \partial_j] = f (\partial_i g) \partial_j - g (\partial_j f) \partial_i.
\end{gather*}
\end{example}

The Schouten bracket $[~,~]_s$ on the shifted complex $C[1] \coloneqq \ChEcochains(\mathfrak{g},\Sym\mathfrak{g})[1]$ induces a bracket $[~,~]_s'$ on the shifted cohomology complex $H[1] \coloneqq \ChE(\mathfrak{g},\Sym\mathfrak{g})[1]$, which is also a graded Lie algebra.

We consider $\overline{d}$ the unique coderivation of $\Sym(H[2])$ induced by $[~,~]'_s[1]$ and $\overline{\delta_{CE} + D}$, the unique coderivation of $\Sym(C[2])$ induced by $\delta_{CE} + [~,~]_s[1]$. They are also differentials, equipping respectively $(H[1],\overline{d})$ and $(C[1],\overline{\delta_{CE}+D})$ of structures of $L_\infty$-algebras.

\begin{defi}
The Chevalley-Eilenberg complex $C[1]$ associated to the Lie algebra $(\mathfrak{g},[~,~])$ is said to be \emph{formal} if there is a $L_\infty$-quasi-isomorphism $\Phi : \Sym(H[2]) \to \Sym(C[2])$, \ie
\begin{equation}
(\Phi \otimes \Phi) \circ \Delta_{\Sym{} H[2]} = \Delta_{\Sym{} C[2]} \circ \Phi \qquad \textrm{et} \qquad \overline{\delta_{CE}+D} \circ \Phi = \Phi \circ \overline{d},
\end{equation}
such that the restriction $\Phi_1$ of $\Phi$ to $C[2]$ is a section. The map $\Phi$ is called a \emph{formality map}.
\end{defi}

Even if there is no formality map between $\Sym(H[2])$ and $\Sym(C[2])$, we can always modify the $L_\infty$ structure induced by $d = d_2 = [~,~]'_s[1]$ in $d = \sum_{n \geqslant 2} d_n$ and by induction build the components of the differential $\overline{d}$ and of the $L_\infty$-morphism $\Phi = \overline{\varphi}$. This construction was done by Bordemann et al. in \cite[A.4, Proposition A.3]{BGHHW05}. This result can also be obtained by a more general point of view using the \emph{Perturbation Lemma} for contraction of differential complexes. This was done by induction by Huebschmann for dg-Lie and $L_\infty$-algebras in \cite{Hueb10,Hueb11}, and Bordemann gives in \cite{B08} a closed formula. \\

The obtained $L_\infty$ morphism $\Phi = \overline{\varphi}$ and new differential $\overline{d}$ satisfy the equation
\begin{equation}
\Phi \circ \overline{d} = \overline{\delta_{CE}+D} \circ \Phi,
\end{equation}
which writes, after projection on $C[2]$ and using that $D = D_2 = [~,~]_s[1]$,
\begin{equation}
\varphi \circ (d \tilde{\star} id_{\Sym(H[2])}) = \delta_{CE} \circ \varphi + \frac{1}{2} D_2 \circ \varphi \tilde{\star} \varphi.
\end{equation}

Since $d = \sum_{n \geqslant 2} d_n$, evaluating on $y_1 \bullet \dotsb \bullet y_{k+1} \in \Sym(H[2])$ gives
\begin{equation} \label{eqs_formality_shifted_twice}
\begin{split}
\sum_{a=2}^{k+1} & \sum_{1 \leqslant i_1 < \dotsb < i_a \leqslant k+1} \nu_a(y_1,\dotsc,y_{k+1}) \\
& \cdotp \varphi_{k+2-a}(d_a(y_{i_1} \bullet \dotsc \bullet y_{i_a}) \bullet y_1 \bullet \dotsc \bullet \widehat{y_{i_1}} \bullet \dotsc \bullet \widehat{y_{i_a}} \bullet \dotsc \bullet y_{k+1}) \\
={} & \delta_{CE} \varphi_{k+1}(y_1 \bullet \dotsc \bullet y_{k+1}) \\
+ \frac{1}{2} & \sum_{a=1}^k \sum_{1 \leqslant i_1 < \dotsb < i_a \leqslant k+1} \nu_a(y_1,\dotsc,y_{k+1}) \\
& \cdotp D_2(\varphi_a(y_{i_1} \bullet \dotsc \bullet y_{i_a}) \bullet \varphi_{k+1-a}(y_1 \bullet \dotsc \bullet \widehat{y_{i_1}} \bullet \dotsc \bullet \widehat{y_{i_a}} \bullet \dotsc \bullet y_{k+1}).
\end{split}
\end{equation}

The maps $\varphi_k$ are of degree $0$, and the signs
\begin{equation*}
\nu_a(y_1,\dotsc,y_{k+1}) = \prod_{r=1}^a (-1)^{|y_{i_r}|(|y_1|+ \dotsb + \widehat{|y_{i_1}|} + \dotsb + \widehat{|y_{i_{r-1}}|} + \dotsb + |y_{i_r-1}|)}
\end{equation*}
for $1 \leqslant a \leqslant k+1$ only come from the permutation which place the elements of indices $i_1,\dotsc,i_a$ at the beginning. \\
However, care of signs must be taken when working with the shifted maps. Denoting by $s : C[1] \to C$ the suspension map which adds one to the degree of elements, we have for example
\begin{equation*}
D_2(x,y) = [x,y]_s[1] = s^{-1} \circ [~,~]_s \circ (s \otimes s)(x \otimes y) = (-1)^{|x|} s^{-1} [s(x),s(y)]_s = (-1)^{|x|}[x,y]_s.
\end{equation*}

\section{Cohomology of the Heisenberg algebra} \label{Sec:Heisenberg-cohomology}

From now on, we take $\mathfrak{g} = \mathfrak{h}_3$ the three-dimensional Heisenberg Lie algebra
\[
\mathfrak{h}_3 = <x,y,z \ \ |\ \  [x,y]=z>
\]
or, writing the Lie bracket as a bivector $[~,~] = \pi = z \partial_x \wedge \partial_y \in \Hom\big({\bigwedge}^2 \mathfrak{g},\mathfrak{g}\big)$.

Since $\dim \mathfrak{g} = 3$, the degree zero part $C^0 = \Sym\mathfrak{g}$ of the Chevalley-Eilenberg complex $\ChEcochains(\mathfrak{g},\Sym\mathfrak{g})$ identifies with $\korps[x,y,z]$. Moreover, the cohomology $\ChE(\mathfrak{g},\Sym\mathfrak{g})$ is concentrated only in degrees $0,1,2,3$. Indeed, $C^k = \Hom\big({\bigwedge}^k \mathfrak{g},\Sym\mathfrak{g}\big) = \{0\}$ for $k \geqslant 4$, because $\partial_x \wedge \partial_y \wedge \partial_z \wedge \partial_t = 0$ for $t \in \{x,y,z\}$. \\

In the following computations, we mostly use the three relations satisfied by the Schouten bracket and the exterior multiplication $\wedge$, namely the graded antisymmetry \eqref{gr-antisymmetry}, the graded Jacobi identity \eqref{gr-Jacobi} and the graded Leibniz rule \eqref{gr-Leibniz}.

\addtocounter{subsection}{-1}

\subsection{\texorpdfstring{Degree $0$}{Degree 0}}

In degree $0$, we have $\ChE[0](\mathfrak{g},\Sym\mathfrak{g}) = \Sym\mathfrak{g}^\mathfrak{g} = \{f \in \Sym\mathfrak{g},\ \delta_{CE}^0(f) = 0\}$. Let $f \in \Sym\mathfrak{g} = \korps[x,y,z]$,
\begin{equation*}
\delta_{CE}^0(f) = [\pi,f]_s = [z \partial_x \wedge \partial_y,f]_s = z \partial_y f \partial_x - z \partial_x f \partial_y,
\end{equation*}
so $\delta_{CE}^0(f) = 0 \Leftrightarrow \left\{
\begin{aligned}
z \partial_y f &= 0 \\
-z \partial_x f &= 0
\end{aligned}
\right. \
\Leftrightarrow f(x,y,z) = \varphi(z)$ is a function of $z$ only. \\[1ex]
Thus $\ChE[0](\mathfrak{g},\Sym\mathfrak{g}) = \korps[z]$.

\subsection{\texorpdfstring{Degree $1$}{Degree 1}}

The $1$-cocycles are $Z^1 = \{X : \mathfrak{g} \to \Sym\mathfrak{g},\ \delta_{CE}^1(X) = 0\}$. We write $X \in C^1 = \Hom(\mathfrak{g},\Sym\mathfrak{g})$ as $X = X_1 \partial_x + X_2 \partial_y + X_3 \partial_z$, with $X_i \in \korps[x,y,z]$.
\begin{align*}
\delta_{CE}^1(X) = [\pi,X]_s &= [z \partial_x \wedge \partial_y,X_1 \partial_x + X_2 \partial_y + X_3 \partial_z]_s \\
&= \left(z \partial_x X_1 + z \partial_y X_2 - X_3\right) \partial_x \wedge \partial_y + z \partial_y X_3\ \partial_x \wedge \partial_z - z \partial_x X_3\ \partial_y \wedge \partial_z,
\end{align*}
\begin{equation*}
\text{so}\quad \delta_{CE}^1(X) = 0 \Leftrightarrow \left\{
\begin{aligned}
z \partial_x X_1 + z \partial_y X_2 - X_3 &= 0 \\
z \partial_y X_3 &= 0 \\
- z \partial_x X_3 &= 0
\end{aligned}
\right..
\end{equation*}
The last two equations imply that $X_3(x,y,z) = \varphi(z)$ is a function of $z$ only. The first equation shows then that $X_3 = \varphi(z) = z \psi(z)$ is a multiple of $z$; and can be rewritten as
\begin{equation}
\partial_x X_1 + \partial_y X_2 = \psi. \label{eqZ1}
\end{equation}
Setting
\(
\left\{
\begin{gathered}
X_1 = - \partial_y g + Y_1 \psi \\
X_2 = \partial_x g + Y_2 \psi
\end{gathered}
\right.,
\)
with $g \in \Sym\mathfrak{g}$, solving \eqref{eqZ1} in $X_1,X_2$ is equivalent to solve in $Y_1,Y_2$ the equation
\begin{equation}
\partial_x Y_1 + \partial_y Y_2 = 1. \label{eqZ11}
\end{equation}
The solutions of \eqref{eqZ11} are of the form \(
\left\{
\begin{aligned}
Y_1 &= a x + K_1(y,z) \\
Y_2 &= (1-a) y + K_2(x,z)
\end{aligned}
\right.,
\)
with $a \in \korps$, and a good choice of the function $g$, gives
\(
\left\{
\begin{aligned}
X_1 &= - \partial_y g + a x \psi \\
X_2 &= \partial_x g + (1-a) y \psi
\end{aligned}
\right..
\)
\\[1ex]

We obtained the following expression for $1$-cocycles
\begin{equation}
X = \left(a x \psi(z) - \partial_y g\right) \partial_x + \big((1-a) y \psi(z) + \partial_x g\big) \partial_y + z \psi(z) \partial_z. \label{1cocycle-expr}
\end{equation}

The $1$-coboundaries are $B^1 = \left\{X : \mathfrak{g} \to \Sym\mathfrak{g},\ \exists\ F \in \Sym\mathfrak{g},\ X = \delta_{CE}^0(F)\right\}$. Since we have $\delta_{CE}^0(F) = z \partial_y F \partial_x - z \partial_x F \partial_y$,\quad $W = \partial_y g \partial_x - \partial_x g \partial_y$ is an element of $B^1$ if and only if $g = z F + \text{cst}$.

Writing $g(x,y,z) = g_0(x,y) + z g_1(x,y,z)$, we have
\[
W = \partial_y g \partial_x - \partial_x g \partial_y = \underbrace{\partial_y g_0(x,y) \partial_x - \partial_x g_0(x,y) \partial_y}_{\in\ Z^1 \setminus B^1} + \delta_{CE}^0(g_1).
\]

Since the $z$ multiple part is a coboundary, we will often use subscript $0$ for the constant part in $z$ and subscript $1$ for the $z$-part. We also use the following notations: $X_g \coloneqq - \partial_y g \partial_x + \partial_x g \partial_y$ for the Lie derivative associated to the function $g$ and $D_a~\coloneqq~a x \partial_x + (1-a) y \partial_y + z \partial_z$. The expression \eqref{1cocycle-expr} rewrites
\begin{equation*}
X = \psi D_a + X_{g_0} - \delta_{CE}^0(g_1),
\end{equation*}
hence the description of the degree one part of the cohomology:
\begin{equation}
\ChE[1](\mathfrak{g},\Sym\mathfrak{g}) = \left\{[X] : \mathfrak{g} \to \Sym\mathfrak{g},\ X = \psi(z) D_a + X_{g_0}\right\},
\end{equation}
with $\psi \in \korps[z]$ and $g_0 \in \korps[x,y]$. (The element $[X]$ is the class of cohomology of $X \in Z^1$.) \\

Let $X = \varphi D_a + X_{g_0}$ and $Y = \psi D_b + X_{h_0}$ two elements of $Z^1$, with $\varphi,\psi \in \korps[z]$, $g_0,h_0 \in \korps[x,y]$ and $a,b \in \korps$. Like said before, we note
\(
\left\{
\begin{aligned}
\varphi(z) &= \varphi_0 + z \varphi_1(z) \\
\psi(z) &= \psi_0 + z \psi_1(z)
\end{aligned}
\right..
\)

We have $D_{1/2} = \dfrac{1}{2} x \partial_x + \dfrac{1}{2} y \partial_y + z \partial_z$, and setting
\begin{flalign*}
& & k(x,y,z) &= \left(\dfrac{1}{2}-a\right)x y \varphi_1(z) &
\tilde{k}(x,y,z) &= \left(\dfrac{1}{2}-b\right)x y \psi_1(z) & &
\end{flalign*}
we get
\begin{align*}
X &= \left(a x \varphi_0 - \partial_y g_0 + \dfrac{1}{2} x z \varphi_1(z)\right) \partial_x + \left((1-a) y \varphi_0 + \partial_x g_0 + \dfrac{1}{2} y z \varphi_1(z)\right) \partial_y + z \varphi(z) \partial_z + \delta_{CE}^0(k) \\
Y &= \left(b x \psi_0 - \partial_y h_0 + \dfrac{1}{2} x z \psi_1(z)\right) \partial_x + \left((1-b) y \psi_0 + \partial_x h_0 + \dfrac{1}{2} y z \psi_1(z)\right) \partial_y + z \psi(z) \partial_z + \delta_{CE}^0(\tilde{k}).
\end{align*}

Thus $X$ and $Y$ are in the same cohomology class if, and only if
\begin{align*}
[X] = [Y] &\Leftrightarrow \left\{
\begin{aligned}
a x \varphi_0 - \partial_y g_0 + \dfrac{1}{2} x z \varphi_1(z) &= b x \psi_0 - \partial_y h_0 + \dfrac{1}{2} x z \psi_1(z) \\
(1-a) y \varphi_0 + \partial_x g_0 + \dfrac{1}{2} y z \varphi_1(z) &= (1-b) y \psi_0 + \partial_x h_0 + \dfrac{1}{2} y z \psi_1(z) \\
z \varphi(z) &= z \psi(z)
\end{aligned}
\right. \\
&\Leftrightarrow \varphi = \psi \qquad \text{and} \qquad g_0 - h_0 = (a-b)x y  \varphi_0,
\end{align*}
and $Y \in [X] \Leftrightarrow \psi = \varphi$ and $h_0 = g_0 - (a-b)x y \varphi_0$. \\

So up to changing the function $g_0 \in \korps[x,y]$, when considering an element $[X] \in \ChE[1](\mathfrak{g},\Sym\mathfrak{g})$, we can always assume that $X = X_{g_0} + \varphi D_a$, with the same $a \in \korps$ fixed once for all.

\subsection{\texorpdfstring{Degree $2$}{Degree 2}}

The $2$-cocycles are $Z^2 = \{\xi : \mathfrak{g} \wedge \mathfrak{g} \to \Sym\mathfrak{g},\ \delta_{CE}^2(\xi) = 0\}$. We write $\xi \in C^2 = \Hom(\mathfrak{g} \wedge \mathfrak{g},\Sym\mathfrak{g})$ as $\xi_{12}\ \partial_x \wedge \partial_y + \xi_{13}\ \partial_x \wedge \partial_z + \xi_{23}\ \partial_y \wedge \partial_z$, with $\xi_{ij} \in \korps[x,y,z]$.
\begin{align*}
\delta_{CE}^2(\xi) = [\pi,\xi]_s &= [z \partial_x \wedge \partial_y,\xi_{12}\ \partial_x \wedge \partial_y + \xi_{13}\ \partial_x \wedge \partial_z + \xi_{23}\ \partial_y \wedge \partial_z]_s \\
&= z \left(\partial_x \xi_{13} + \partial_y \xi_{23}\right) \partial_x \wedge \partial_y \wedge \partial_z.
\end{align*}
\begin{equation*}
\text{so}\quad \delta_{CE}^2(\xi) = 0 \Leftrightarrow z \left(\partial_x \xi_{13} + \partial_y \xi_{23}\right) = 0 \Leftrightarrow \xi_{13} = - \partial_y g\quad \text{and}\quad \xi_{23} = \partial_x g
\end{equation*}
We can express the $2$-cocycles as
\begin{equation}
\xi = h\ \partial_x \wedge \partial_y - \partial_y g\ \partial_x \wedge \partial_z + \partial_x g\ \partial_y \wedge \partial_z,
\end{equation}
with $g,h$ arbitrary functions in $\korps[x,y,z]$. \\

The $2$-coboundaries are $B^2 = \left\{\xi : \mathfrak{g} \wedge \mathfrak{g} \to \Sym\mathfrak{g},\ \exists\ X \in \ChEcochains[1](\mathfrak{g},\Sym\mathfrak{g}),\ \xi = \delta_{CE}^1(X)\right\}$.
Since $\delta_{CE}^1(X) =  \left(z \partial_x X_1 + z \partial_y X_2 - X_3\right) \partial_x \wedge \partial_y + z \partial_y X_3\ \partial_x \wedge \partial_z - z \partial_x X_3\ \partial_y \wedge \partial_z,$
\[
\xi = \delta_{CE}^1(X) \Leftrightarrow \left\{
\begin{aligned}
h &= z \partial_x X_1 + z \partial_y X_2 - X_3 \\
- \partial_y g &= z \partial_y X_3 \\
\partial_x g &= z \partial_x X_3
\end{aligned}
\right.
\]

We note
\begin{equation*}
\begin{aligned}
g(x,y,z) &= g_0(x,y) + z g_1(x,y,z) \\
&= g_0(x,y) + z \left(g_{10}(x,y) + z g_{11}(x,y,z)\right)
\end{aligned} \qquad \text{and} \qquad
\begin{aligned}
h(x,y,z) &= \widetilde{h_0}(x,y) + z h_1(x,y,z) \\
h_0 &= \widetilde{h_0} - g_{10},
\end{aligned}
\end{equation*}
and we take $X = X_1 \partial_x + X_2 \partial_y + X_3 \partial_z$ with $X_3 = g_1$ and $X_1,X_2$ such that $\partial_x X_1 + \partial_y X_2 = g_{11} - h_1$.
We then have
\[
\xi + \delta_{CE}^1(X) = h_0\ \partial_x \wedge \partial_y - \partial_y g_0\ \partial_x \wedge \partial_z + \partial_x g_0\ \partial_y \wedge \partial_z,
\]
with $g_0,h_0$ function of $x$ and $y$ only. Further setting $Y = h_0 \partial_z$, we get
\begin{align*}
\xi + \delta_{CE}^1(X) + \delta_{CE}^1(Y) &= - \partial_y (g_0 - z h_0)\ \partial_x \wedge \partial_z + \partial_x (g_0 - z h_0)\ \partial_y \wedge \partial_z \\
&= X_{g_0 - z h_0} \wedge \partial_z,
\end{align*}
recalling that $X_f = - \partial_y f \partial_x + \partial_x f \partial_y$.

Hence the description of the degree two part of the cohomology:
\begin{equation}
\ChE[2](\mathfrak{g},\Sym\mathfrak{g}) = \left\{[\xi] : \mathfrak{g} \wedge \mathfrak{g} \to \Sym\mathfrak{g},\ \xi = X_g \wedge \partial_z\right\},
\end{equation}
with $g \in \korps[x,y,z]$.

\subsection{\texorpdfstring{Degree $3$}{Degree 3}}

Since the Schouten bracket $[~,~]_s$ is of degree $0$, for $\xi \in C^k$, $\delta_{CE}(\xi) = [\pi,\xi]_s$ is of degree $1+k-1=k$. (Recall that in $C[1]$, elements of $C^k = \ChEcochains[k](\mathfrak{g},\Sym\mathfrak{g})$ are of degree $k-1$.)
Therefore, all elements of $C^3=\Hom\big({\bigwedge}^3 \mathfrak{g},\Sym\mathfrak{g}\big)$ are cocycles, since they are of degree $3$ and that $C^4 = \{0\}$, thus $Z^3 = \ChEcochains[3](\mathfrak{g},\Sym\mathfrak{g})$. \\

The $3$-coboundaries are $B^3 = \left\{\xi : \mathfrak{g} \wedge \mathfrak{g} \wedge \mathfrak{g} \to \Sym\mathfrak{g},\ \exists\ \phi \in \ChEcochains[2](\mathfrak{g},\Sym\mathfrak{g}),\ \xi = \delta_{CE}^2(\phi)\right\}$.
\begin{gather*}
\xi = \delta_{CE}^2(\phi) \Leftrightarrow \xi_{123}\ \partial_x \wedge \partial_y \wedge \partial_z = z \left(\partial_x \phi_{13} + \partial_y \phi_{23}\right) \partial_x \wedge \partial_y \wedge \partial_z
\end{gather*}
Setting $\xi_{123} = \xi_0(x,y) + z \xi_1(x,y,z)$, we take $\phi = \phi_{12}\ \partial_x \wedge \partial_y + \phi_{13}\ \partial_x \wedge \partial_z + \phi_{23}\ \partial_y \wedge \partial_z$ such that $\partial_x \phi_{13} + \partial_y \phi_{23} = - \xi_1$, so that
\begin{equation*}
\xi + \delta_{CE}^2(\phi) = \xi_0\ \partial_x \wedge \partial_y \wedge \partial_z.
\end{equation*}

Hence the description of the degree three part of the cohomology:
\begin{equation}
\ChE[3](\mathfrak{g},\Sym\mathfrak{g}) = \left\{[\xi] : \mathfrak{g} \wedge \mathfrak{g} \wedge \mathfrak{g} \to \Sym\mathfrak{g},\ \xi = \xi_0\ \partial_x \wedge \partial_y \wedge \partial_z\right\},
\end{equation}
with $\xi_0 \in \korps[x,y]$.

\section{Schouten brackets on cochain level} \label{Schouten-cochain}

We work with the graded Lie algebra $C[1] = \ChEcochains(\mathfrak{g},\Sym\mathfrak{g})[1]$ where $\mathfrak{g} = \mathfrak{h}_3$. We compute in this section the various Schouten brackets $[C^i,C^j]_s$, for $0 \leqslant i \leqslant j \leqslant 3$. In the computations, we do not use the Lie structure of the Heisenberg algebra, all that is used is that $\dim \mathfrak{g} = 3$. \\

Reasoning on degrees, for $f,g \in C^0 = \Sym\mathfrak{g}$, we have
\begin{equation} \label{cochain00}
[f,g]_s = 0.
\end{equation}

\medskip

Let $f \in C^0$ and $X = X_1 \partial_x + X_2 \partial_y + X_3 \partial_z \in C^1$, with $X_i \in \korps[x,y,z]$.
\begin{equation} \label{cochain01}
[X,f]_s = X(f) = X_1 \partial_x f + X_2 \partial_y f + X_3 \partial_z f
\end{equation}

\medskip

Let $f \in C^0$ and $\xi = \xi_{12}\ \partial_x \wedge \partial_y + \xi_{13}\ \partial_x \wedge \partial_z + \xi_{23}\ \partial_y \wedge \partial_z \in C^2$, with $\xi_{ij} \in \korps[x,y,z]$.
\begin{equation} \label{cochain02}
[f,\xi]_s = \xi_{12} \left(\partial_y f \partial_x - \partial_x f \partial_y\right) + \xi_{13} \left(\partial_z f \partial_x - \partial_x f \partial_z\right) + \xi_{23} \left(\partial_z f \partial_y - \partial_y f \partial_z\right)
\end{equation}

\medskip

Let $f \in C^0$ and $\xi = \xi_{123}\ \partial_x \wedge \partial_y \wedge \partial_z \in C^3$, with $\xi_{123} \in \korps[x,y,z]$.
\begin{equation} \label{cochain03}
[f,\xi]_s = \xi_{123} \left(- \partial_z f\ \partial_x \wedge \partial_y + \partial_y f\ \partial_x \wedge \partial_z - \partial_x f\ \partial_y \wedge \partial_z\right)
\end{equation}

\medskip

Let $X,Y \in C^1$, with $X = X_1 \partial_x + X_2 \partial_y + X_3 \partial_z$ and $Y = Y_1 \partial_x + Y_2 \partial_y + Y_3 \partial_z$
\begin{align} \label{cochain11}
[X,Y]_s ={}& \quad\!\! \left(X_1 \partial_x Y_1 - Y_1 \partial_x X_1 + X_2 \partial_y Y_1 - Y_2 \partial_y X_1 + X_3 \partial_z Y_1 - Y_3 \partial_z X_1\right) \partial_x \notag \\
&+ \left(X_1 \partial_x Y_2 - Y_1 \partial_x X_2 + X_2 \partial_y Y_2 - Y_2 \partial_y X_2 + X_3 \partial_z Y_2 - Y_3 \partial_z X_2\right) \partial_y \\
&+ \left(X_1 \partial_x Y_3 - Y_1 \partial_x X_3 + X_2 \partial_y Y_3 - Y_2 \partial_y X_3 + X_3 \partial_z Y_3 - Y_3 \partial_z X_3\right) \partial_y \notag
\end{align}

\medskip

Let $X = X_1 \partial_x + X_2 \partial_y + X_3 \partial_z \in C^1$ and $\xi = \xi_{12}\ \partial_x \wedge \partial_y + \xi_{13}\ \partial_x \wedge \partial_z + \xi_{23}\ \partial_y \wedge \partial_z \in C^2$.
\begin{adjustwidth}{-1cm}{-1cm} \vspace*{-1em}
\begin{align} \label{cochain12}
[X,\xi]_s ={}& \quad\!\! \left(X_1 \partial_x \xi_{12} - \xi_{12} \partial_x X_1 + X_2 \partial_y \xi_{12} - \xi_{12} \partial_y X_2 + X_3 \partial_z \xi_{12} - \xi_{13} \partial_z X_2 + \xi_{23} \partial_z X_1\right) \partial_x \wedge \partial_y \notag \\
&+ \left(X_1 \partial_x \xi_{13} - \xi_{13} \partial_x X_1 + X_2 \partial_y \xi_{13} - \xi_{12} \partial_y X_3 + X_3 \partial_z \xi_{13} - \xi_{13} \partial_z X_3 - \xi_{23} \partial_y X_1\right) \partial_x \wedge \partial_z \\
&+ \left(X_1 \partial_x \xi_{23} - \xi_{13} \partial_x X_2 + X_2 \partial_y \xi_{23} - \xi_{23} \partial_y X_2 + X_3 \partial_z \xi_{23} - \xi_{23} \partial_z X_3 + \xi_{12} \partial_x X_3\right) \partial_y \wedge \partial_z \notag
\end{align}
\end{adjustwidth}

\medskip

Let $X = X_1 \partial_x + X_2 \partial_y + X_3 \partial_z \in C^1$ and $\xi = \xi_{123}\ \partial_x \wedge \partial_y \wedge \partial_z \in C^3$.
\begin{equation} \label{cochain13}
[X,\xi]_s = \big(X_1 \partial_x \xi_{123} + X_2 \partial_y \xi_{123} + X_3 \partial_z \xi_{123} -  \xi_{123} \left(\partial_x X_1 + \partial_y X_2 + \partial_z X_3\right) \big) \partial_x \wedge \partial_y \wedge \partial_z
\end{equation}

\medskip

Let $\xi,\eta \in C^2$, with $\xi = \xi_{12}\ \partial_x \wedge \partial_y + \xi_{13}\ \partial_x \wedge \partial_z + \xi_{23}\ \partial_y \wedge \partial_z$ and $\eta = \eta_{12}\ \partial_x \wedge \partial_y + \eta_{13}\ \partial_x \wedge \partial_z + \eta_{23}\ \partial_y \wedge \partial_z$.
\begin{equation} \label{cochain22}
\begin{split}
[\xi,\eta]_s ={}& \bigg( \xi_{12} \partial_x \eta_{13} - \eta_{13} \partial_x \xi_{12} + \xi_{12} \partial_y \eta_{23} - \eta_{23} \partial_y \xi_{12} + \eta_{12} \partial_x \xi_{13} - \xi_{13} \partial_x \eta_{12} \\
\quad & + \xi_{13} \partial_z \eta_{23} - \eta_{23} \partial_z \xi_{13} + \eta_{12} \partial_y \xi_{23} - \xi_{23} \partial_y \eta_{12} + \eta_{13} \partial_z \xi_{23} - \xi_{23} \partial_z \eta_{13} \bigg) \partial_x \wedge \partial_y \wedge \partial_z
\end{split}
\end{equation}
Reasoning on degrees, we already had that $[\xi,\eta]_s = 0$ for $\xi \in C^2$ and $\eta \in C^3$ and it is also the case for $\xi \in C^3$ and $\eta \in C^3$.

\section{Schouten brackets on cohomology} \label{Schouten-cohomology}

We will now compute the various Schouten brackets $[C^i,C^j]'_s$ on cohomology, for $0 \leqslant i \leqslant j \leqslant 3$. Since we use cocycle representatives, we will actually do the computations with the Schouten bracket $[~,~]_s$ on cochain level by taking care of coboundaries. When writing elements of the cohomology, we will abusively omit the class notation to lighten the expressions. The equations on cochain level with the coboundaries will be useful for expressing formality equations in \autoref{Sec:formality-eq}. \\

Let $\varphi,\psi \in H^0 = \ChE[0](\mathfrak{g},\Sym\mathfrak{g}) = \korps[z]$. Like for cochains,
\begin{equation} \label{cohomology00}
[\varphi,\psi]'_s = 0.
\end{equation}

\medskip

Let $\varphi \in H^0$ and $X_{g_0} + \psi D_a \in H^1$, with $\psi \in \korps[z]$ and $g_0 \in \korps[x,y]$.
\begin{equation} \label{cohomology01}
[\varphi,X_{g_0} + \psi D_a]_s = [\varphi,z \psi \partial_z]_s = - z \psi \partial_z \varphi
\end{equation}

\medskip

Let $\varphi \in H^0$ and $X_g \wedge \partial_z \in H^2$, with $g \in \korps[x,y,z]$.
\begin{equation} \label{cohomology02}
[\varphi,X_g \wedge \partial_z]_s = X_{g \varphi'_z} = X_{(g \varphi'_z)_0} + \delta_{CE}^0\left((g \varphi'_z)_1\right)
\end{equation}
Here, $\varphi'_z$ is the derivative of $\varphi$ with respect to $z$, and $g \varphi'_z = (g \varphi'_z)_0 + z (g \varphi'_z)_1$ are the constant term (in $z$) and remaining $z$-part. Hence on the cohomology, we have
\begin{equation}
[\varphi,X_g \wedge \partial_z]'_s = X_{(g \varphi'_z)_0} = X_{g_0 \varphi'_z(0)}.
\end{equation}

\medskip

Let $\varphi \in H^0$ and $P\ \partial_x \wedge \partial_y \wedge \partial_z \in H^3$, with $P \in \korps[x,y]$.
\begin{equation} \label{cohomology03}
\begin{split}
[\varphi,P\ \partial_x \wedge \partial_y \wedge \partial_z]_s ={} & - P \varphi'_z \partial_x \wedge \partial_y \\
={} & \delta_{CE}^1(P \varphi'_z \partial_z) + X_{P z \varphi'_z} \wedge \partial_z
\end{split}
\end{equation}
Hence on the cohomology,
\begin{equation}
[\varphi,P\ \partial_x \wedge \partial_y \wedge \partial_z]'_s = X_{P z \varphi'_z} \wedge \partial_z.
\end{equation}

\medskip

Let $X_{g_0} + \varphi D_a,\ X_{h_0} + \psi D_a \in H^1$, with $g_0,h_0 \in \korps[x,y]$ and $\varphi,\psi \in \korps[z]$.
\begin{equation} \label{cohomology11}
\begin{split}
[X_{g_0} &+ \varphi D_a, X_{h_0} + \psi D_a]_s \\ ={}& X_{\{g_0,h_0\} + \varphi_0 (D_a(h_0) - h_0) - \psi_0 (D_a(g_0) - g_0)} + z(\varphi \psi'_z - \varphi'_z \psi) D_a \\
&+ \delta_{CE}^0\big(\psi'_z (D_a(g_0) - g_0) - \varphi'_z (D_a(h_0) - h_0)\big)
\end{split}
\end{equation}
The Poisson bracket $\{~,~\}$ is the one for functions of two variables $x$ and $y$: we note $\{g,h\} = \partial_x g \partial_y h - \partial_x h \partial_y g$; and we have $[X_g,X_h] = X_{\{g,h\}}$.
On cohomology we get
\begin{equation}
\begin{split}
[X_{g_0} &+ \varphi D_a, X_{h_0} + \psi D_a]'_s \\ ={}& X_{\{g_0,h_0\} + \varphi_0 (D_a(h_0) - h_0) - \psi_0 (D_a(g_0) - g_0)} + z(\varphi \psi'_z - \varphi'_z \psi) D_a \\
\eqqcolon & X_{g_0,\varphi,h_0,\psi} + z(\varphi \psi'_z - \varphi'_z \psi) D_a,
\end{split}
\end{equation}
with the notation $(g_0,\varphi,h_0,\psi)$ standing for the long function in subscript.

\medskip

Let $X_{g_0} + \varphi D_a \in H^1$ and $X_p \wedge \partial_z \in H^2$, with $g_0 \in \korps[x,y]$, $\varphi \in \korps[z]$ and $p \in \korps[x,y,z]$.
\begin{equation} \label{cohomology12}
\begin{split}
[X_{g_0} + \varphi D_a,X_p \wedge \partial_z]_s ={}& X_{\{g_0,p\} + \varphi (D_a(p) - p) - p (\varphi'_z + \varphi) - z \varphi'_z (D_a(p) - z \partial_z p)} \wedge \partial_z \\
& + \delta_{CE}^1\big(- \varphi'_z (D_a(p) - z \partial_z p) \partial_z\big)
\end{split}
\end{equation}
So we get
\begin{equation}
\begin{split}
[X_{g_0} + \varphi D_a,X_p \wedge \partial_z]_s ={}& X_{\{g_0,p\} + \varphi (D_a(p) - p) - p (\varphi'_z + \varphi) - z \varphi'_z (D_a(p) - z \partial_z p)} \wedge \partial_z \\
\eqqcolon & X_{g_0,\varphi,p} \wedge \partial_z
\end{split}
\end{equation}
with the notation $(g_0,\varphi,p)$ standing for the long function in subscript. (It cannot be confused with the preceding notation since it has only three arguments.)

\medskip

Let $X_{g_0} + \varphi D_a \in H^1$ and $P\ \partial_x \wedge \partial_y \wedge \partial_z \in H^3$, with $g_0,P \in \korps[x,y]$ and $\varphi \in \korps[z]$.
\begin{equation} \label{cohomology13}
\begin{split}
[X_{g_0} + \varphi D_a,P\ \partial_x \wedge \partial_y \wedge \partial_z]_s ={}& \left(X_{g_0}(P) + \varphi (D_a(P) - 2P) - P z \varphi'_z\right) \partial_x \wedge \partial_y \wedge \partial_z \\
={}& \left(X_{g_0}(P) + \varphi_0 (D_a(P) - 2P)\right) \partial_x \wedge \partial_y \wedge \partial_z \\
& + \delta_{CE}^2(\phi),
\end{split}
\end{equation}
with $\phi = \phi_{13}\ \partial_x \wedge \partial_z + \phi_{23} \partial_y \wedge \partial_z$ such that $\partial_x \phi_{13} + \partial_y \phi_{23} = \varphi_1 (D_a(P) - 2P) - P \varphi'_z$.
So we have on the cohomology
\begin{equation}[X_{g_0} + \varphi D_a,P\ \partial_x \wedge \partial_y \wedge \partial_z]'_s = \left(X_{g_0}(P) + \varphi_0 (D_a(P) - 2P)\right) \partial_x \wedge \partial_y \wedge \partial_z.
\end{equation}

\medskip

Let $X_f \wedge \partial_z,\ X_g \wedge \partial_z \in H^2$, with $f,g \in \korps[x,y,z]$. Writing
\(
\left\{
\begin{aligned}
f(x,y,z) &= f_0(x,y) + z f_1(x,y,z) \\
g(x,y,z) &= g_0(x,y) + z g_1(x,y,z)
\end{aligned}
\right.
\) we get
\begin{equation} \label{cohomology22}
[X_f \wedge \partial_z,X_g \wedge \partial_z]_s = [X_f \wedge \partial_z,X_g \wedge \partial_z]'_s = \big(\{f_0,g_1\} + \{g_0,f_1\}\big)\partial_x \wedge \partial_y \wedge \partial_z.
\end{equation}

For completeness, we write the last two cases were the brackets are trivial. Let $X_f \wedge \partial_z \in H^2$ and $P\ \partial_x \wedge \partial_y \wedge \partial_z,\ Q\ \partial_x \wedge \partial_y \wedge \partial_z \in H^3$, with $f \in \korps[x,y,z]$ and $P,Q \in \korps[x,y]$.
\begin{align}\label{cohomology233}
[X_f \wedge \partial_z,P\ \partial_x \wedge \partial_y \wedge \partial_z]'_s ={}& 0 \\
[P\ \partial_x \wedge \partial_y \wedge \partial_z,Q\ \partial_x \wedge \partial_y \wedge \partial_z]'_s ={}& 0
\end{align}

\section{Formality equations} \label{Sec:formality-eq}

In this section, we will explicit the formality equations \eqref{eqs_formality_shifted_twice} for the first orders. In particular, we will obtain non-trivial terms $d_k$ of arity $k=3$ and higher showing that the Chevalley-Eilenberg complex of the Heisenberg algebra $\mathfrak{h}_3$ is not formal in Kontsevich sense, \ie the $L_\infty$ structure on the cohomology is not induced by the Schouten bracket $d_2 = [~,~]'_s[1]$ only.

Note that in equations \eqref{eqs_formality_shifted_twice}, we work with elements $y_i \in H[2]$ so the degrees are shifted again. \\

The maps $\varphi_k : \Sym^k(H[2]) \to C[2]$ are of degree zero, so for $i_1,\dotsc,i_n \in \text{\textlbrackdbl}-2,1\text{\textrbrackdbl}$
\begin{equation*}
\varphi_k(H[2]^{i_1} \bullet \dotsb \bullet H[2]^{i_n}) \subset C[2]^{i_1+\dotsb+i_n}.
\end{equation*}

Sometimes, we will consider the maps $\varphi_k : H[2]^{\otimes k} \to C[2]$ having source in the tensorial power instead of the symmetric algebra of the space, and write $\varphi_k(H[2]^{i_1},\dotsc,H[2]^{i_n})$ with periods to separate arguments.

Since the cohomology is non-trivial in degrees $0,1,2,3$, we can always choose $k$ elements in $H[2]$ of various degrees such that the image by the map $\varphi_k$ is in the non-trivial part of $C[2]$, it suffice to choose $k$ numbers in $\text{\textlbrackdbl}-2,1\text{\textrbrackdbl}$ such that their sum is again in $\text{\textlbrackdbl}-2,1\text{\textrbrackdbl}$. The situation here is more complicated that for the Lie algebra $\mathfrak{so}(3)$ (see \cite{El12}), since in this case we have terms $\varphi_k$, and hence terms $d_k$ of all arities greater than two.

\addtocounter{subsection}{-1}

\subsection{\texorpdfstring{Order $k=0$}{Order k=0}}
For $k=0$ and $y_1 \in H[2]$, we have
\begin{equation}
0 = \delta_{CE}(\varphi_1(y_1))
\end{equation}
which express that $\varphi_1$ is a section from the cohomology into the cocycles. We choose $\varphi_1 = \operatorname{inc}$, the inclusion given by the choice of representatives of classes given in \autoref{Sec:Heisenberg-cohomology}.

\subsection{\texorpdfstring{Order $k=1$}{Order k=1}}
For $k=1$ and $y_1 \bullet y_2 \in \Sym^2(H[2])$, we have
\begin{equation} \label{formality-eq-k1}
\varphi_1(d_2(y_1 \bullet y_2)) = \delta_{CE} \varphi_2(y_1 \bullet y_2) + D(\varphi_1(y_1) \bullet \varphi_1(y_2)).
\end{equation}
This equation shows that the inclusion $\varphi_1$ is not a morphism of the graded Lie algebras $(H[2],d_2) \to (C[2],D)$, since there is the coboundary $\delta_{CE} \varphi_2(y_1 \bullet y_2)$. We will explicit the term $\varphi_2$ and higher terms to obtain an $L_\infty$-morphism between these two $L_\infty$-algebras.
To simplify the expressions, we will omit the inclusion $\varphi_1$ in the remaining computations. The preceding equation rewrites
\begin{equation}
d_2(y_1 \bullet y_2) = \delta_{CE} \varphi_2(y_1 \bullet y_2) + D(y_1 \bullet y_2).
\end{equation} \\

Let $y_1=\varphi,y_2=\psi \in H^0 = \korps[z]$. We have
\(
d_2(\varphi \bullet \psi) = (-1)^{-2} [\varphi,\psi]'_s = 0
\)
and also
\(
D(\varphi \bullet \psi) = (-1)^{-2} [\varphi,\psi]'_s = 0
\)
so $\delta_{CE} \varphi_2(\varphi,\psi) = 0$ and we can set 
\begin{equation}
\varphi_2(\varphi,\psi) = 0.
\end{equation}

\medskip

Let $y_1 = \varphi \in H^0$ and $y_2 = X_{g_0} + \psi D_a \in H^1$ with $\varphi,\psi \in \korps[z]$ and $g_0 \in \korps[x,y]$. As before
\(
[\varphi,X_{g_0} + \psi D_a]_s = [\varphi,X_{g_0} + \psi D_a]'_s = - z \psi \varphi'_z,
\)
and we set again
\begin{equation}
\varphi_2(\varphi,X_{g_0} + \psi D_a) = 0.
\end{equation}

\medskip

Let $y_1 = \varphi \in H^0$ and $y_2 = X_g \wedge \partial_z \in H^2$, with $\varphi \in \korps[z]$ and $g \in \korps[x,y,z]$. Equation \eqref{cohomology02} rewrites
\(
D(\varphi \bullet X_g \wedge \partial_z) = d_2(\varphi \bullet X_g \wedge \partial_z) + \delta_{CE}^0((g\varphi'_z)_1),
\)
so to obtain the formality equation \eqref{formality-eq-k1} in this case, we can set
\begin{equation}
\varphi_2(\varphi,X_g \wedge \partial_z) = - (g\varphi'_z)_1.
\end{equation}

\medskip

Let $y_1 = \varphi \in H^0$ and $y_3 = P \partial_x \wedge \partial_y \wedge \partial_z \in H^3$, with $\varphi \in \korps[z]$ and $P \in \korps[x,y]$. Equation \eqref{cohomology03} rewrites
\(
D(\varphi \bullet P \partial_x \wedge \partial_y \wedge \partial_z) = d_2(\varphi \bullet P \partial_x \wedge \partial_y \wedge \partial_z) + \delta_{CE}^1(P \varphi'_z \partial_z),
\)
so we set 
\begin{equation}
\varphi_2(\varphi,P \partial_x \wedge \partial_y \wedge \partial_z) = - P \varphi'_z \partial_z.
\end{equation}

\medskip

Let $y_1 = X_{g_0} + \varphi D_a,\ y_2 = X_{h_0} + \psi D_a \in H^1$, with $g_0,h_0 \in \korps[x,y]$ and $\varphi,\psi \in \korps[z]$. As in the previous cases, equation \eqref{cohomology11} gives
\begin{equation}
\varphi_2(X_{g_0} + \varphi D_a,X_{h_0} + \psi D_a) = \psi'_z (D_a(g_0) - g_0) - \varphi'_z (D_a(h_0) - h_0).
\end{equation}

\medskip

For $y_1 = X_{g_0} + \varphi D_a \in H^1$ and $y_2 = X_p \wedge \partial_z \in H^2$, with $\varphi \in \korps[z]$ and $g_0,p \in \korps[x,y]$, equation \eqref{cohomology12} gives
\begin{equation}
\varphi_2(X_{g_0} + \varphi D_a,X_p \wedge \partial_z) = - \varphi'_z (D_a(p) - z \partial_z p) \partial_z
\end{equation}

\medskip

For $y_1 = X_{g_0} + \varphi D_a \in H^1$ and $y_2 = P\ \partial_x \wedge \partial_y \wedge \partial_z \in H^3$, with $\varphi \in \korps[z]$ and $g_0,P \in \korps[x,y]$, equation \eqref{cohomology13} gives 
\begin{equation}
\varphi_2(X_{g_0} + \varphi D_a,P\ \partial_x \wedge \partial_y \wedge \partial_z) = \phi_{13}\ \partial_x \wedge \partial_z + \phi_{23} \partial_y \wedge \partial_z
\end{equation}
such that $\partial_x \phi_{13} + \partial_y \phi_{23} = \varphi_1 (D_a(P) - 2P) - P \varphi'_z$.

\medskip

Let $X_f \wedge \partial_z,\ X_g \wedge \partial_z \in H^2$, with $f,g \in \korps[x,y,z]$, and $P\ \partial_x \wedge \partial_y \wedge \partial_z,\ Q\ \partial_x \wedge \partial_y \wedge \partial_z \in H^3$, with $P,Q \in \korps[x,y]$. Since in equation \eqref{cohomology22} the brackets are the same on cochain level and cohomology, we have
\begin{equation}
\varphi_2(X_f \wedge \partial_z,X_g \wedge \partial_z) = 0,
\end{equation}
and the last two brackets are trivial in both levels, so
\begin{align}
\varphi_2(X_f \wedge \partial_z,P\ \partial_x \wedge \partial_y \wedge \partial_z) = 0 \\
\varphi_2(P\ \partial_x \wedge \partial_y \wedge \partial_z,Q\ \partial_x \wedge \partial_y \wedge \partial_z) = 0
\end{align}

\subsection{\texorpdfstring{Order $k=2$}{Order k=2}}
For $k=2$ and $y_1 \bullet y_2 \bullet y_3 \in \Sym^3(H[2])$, we have
\begin{equation} \label{formality-eq-k2}
\begin{split}
d_3(& y_1 \bullet y_2 \bullet y_3) \\
& + \varphi_2(d_2(y_1 \bullet y_2) \bullet y_3) + (-1)^{|y_2| |y_3|} \varphi_2(d_2(y_1 \bullet y_3) \bullet y_2) + (-1)^{|y_1|(|y_2|+|y_3|)} \varphi_2(d_2(y_2 \bullet y_3) \bullet y_1) \\
={} & \delta_{CE} \varphi_3(y_1 \bullet y_2 \bullet y_3) \\
& + D(y_1 \bullet \varphi_2(y_2 \bullet y_3)) + (-1)^{|y_1| |y_2|} D(y_2 \bullet \varphi_2(y_1 \bullet y_3)) + (-1)^{|y_3|(|y_1|+|y_2|)} D(y_3 \bullet \varphi_2(y_1 \bullet y_2)).
\end{split}
\end{equation}

Recall that the maps $\delta_{CE}$, $D$ and $d_k$ are of degree~$1$ and the maps $\varphi_k$ are of degree~$0$. Thus, the left-hand side and the right-hand side of the equation are of degree one more than the sum of the degrees of the three elements $y_1,y_2,y_3$.

We will keep the same notation as before for elements of various degrees in $H[2]$, without specifying further the types of the functions used (of $x$ and $y$, of $z$ only, of $x,y,z$).
Although we work with double-shifted spaces, we keep the original degree for the Chevalley-Eilenberg coboundary operator, so
\[
\delta_{CE}^k : C^k = \ChEcochains[k](\mathfrak{g},\Sym\mathfrak{g}) \to C^{k+1} = \ChEcochains[k+1](\mathfrak{g},\Sym\mathfrak{g}).
\]
By graded symmetry, it will be enough to study the cases $-2 \leqslant i \leqslant j \leqslant k \leqslant 1$, where $i,j,k$ are the degrees of $y_1,y_2,y_3 \in H[2]$. We note $\omega \coloneqq \partial_x \wedge \partial_y \wedge \partial_z$. \\

Let $y_1 = \varphi, y_2 = \psi \in H^0$ and $P \omega \in H^3$. The equation \eqref{formality-eq-k2} rewrites
\begin{equation*}
d_3(\varphi \bullet \psi \bullet P \omega) - 2(Pz\varphi'_z\psi'_z)_1 = \delta_{CE} \varphi_3(\varphi \bullet \psi \bullet P \omega) + 2P\varphi'_z\psi'_z
\end{equation*}
The elements of this equation are of degree $|\varphi|+|\psi|+|P \omega| + 1 = -2-2+1+1 = -2$ and since $\delta_{CE}$ is of degree $+1$, the element $\varphi_3(\varphi \bullet \psi \bullet P \omega) \in C[2]^{-3} = C^{-1}$ has image in $C^0$, so here $\delta_{CE} = \delta_{CE}^{-1} : C^{-1} \to C^0$. Of course, $C^{-1} = \{0\}$, so we get
\begin{equation*}
\varphi_3(\varphi \bullet \psi \bullet P \omega) = 0 \quad \text{and} \quad d_3(\varphi \bullet \psi \bullet P \omega) = 4P\varphi'_z\psi'_z.
\end{equation*}
Even if we had chosen other expressions for $\varphi_2(\varphi,P \omega)$ and $\varphi_2(\varphi,X_g \wedge \partial_z)$ adding cocycles, there is no possibility to cancel the constant part in $z$ coming from the function $P$. Therefore, this case alone shows that the Chevalley-Eilenberg complex of the Heisenberg algebra is not formal. \\

Let $y_1 = \varphi \in H^0$, $y_2 = X_{g_0} + \psi D_a \in H^1$ and $y_3 = X_h \wedge \partial_z \in H^2$. The equation \eqref{formality-eq-k2} rewrites
\begin{gather*}
\begin{aligned}
d_3(& \varphi \bullet X_{g_0} + \psi D_a \bullet X_h \wedge \partial_z) + \psi'_z \big(D_a(h_0 \varphi'_z(0))-h_0 \varphi'_z(0)\big) + \big((g_0,\psi,h)\varphi'_z\big)_1 \\
={}& \delta_{CE}^{-1}(\varphi \bullet X_{g_0} + \psi D_a \bullet X_h \wedge \partial_z) - \psi'_z \varphi'_z(z\partial_z h-D_a(h))-z \psi \partial_z (h \psi'_z)_1
\end{aligned} \\[1ex]
\begin{aligned}
\!\!\Leftrightarrow d_3(\varphi \bullet X_{g_0} + \psi D_a \bullet X_h \wedge \partial_z) ={}& \psi'_z \big(h_0 \varphi'_z(0) - D_a(h_0 \varphi'_z(0)) - \varphi'_z(z\partial_z h-D_a(h))\big) \\
& - \big((g_0,\psi,h)\varphi'_z\big)_1 - z \psi \partial_z (h \psi'_z)_1.
\end{aligned}
\end{gather*}

Let $y_1 = \varphi \in H^0$, $y_2 = X_{g_0} + \psi D_a \in H^1$ and $y_3 = P \omega \in H^3$. The equation \eqref{formality-eq-k2} rewrites
\begin{equation*}
\begin{split}
d_3(&\varphi \bullet X_{g_0} + \psi D_a \bullet P \omega) \\ ={}& \delta_{CE}^0 \varphi_3(\varphi \bullet X_{g_0} + \psi D_a \bullet P \omega) - \varphi'_z (P \psi'_z a x + \phi_{13}) \partial_x - \varphi'_z (P \psi'_z (1-a) y + \phi_{23}) \partial_y \\
& + \big( X_{g_0}(P)(1-\varphi'_z) + \left(\psi - \psi_0 \varphi'_z - \psi'_z \varphi'_z z\right) D_a(P) - 2 P \varphi'_z(\psi + \psi_0) - \psi'_z \varphi'_z z P \big) \partial_z,
\end{split}
\end{equation*}
where $\phi_{13},\phi_{23}$ come from
\begin{equation*}
\varphi_2(X_{g_0} + \psi D_a,P\ \partial_x \wedge \partial_y \wedge \partial_z) = \phi_{13}\ \partial_x \wedge \partial_z + \phi_{23} \partial_y \wedge \partial_z
\end{equation*}
with $\partial_x \phi_{13} + \partial_y \phi_{23} = \psi_1 (D_a(P) - 2P) - P \psi'_z$.

We can choose the value of $\varphi_3(\varphi \bullet X_{g_0} + \psi D_a \bullet P \omega)$ such that its image by $\delta_{CE}^0$ cancels the coboundary $- \varphi'_z (P \psi'_z a x + \phi_{13}) \partial_x - \varphi'_z (P \psi'_z (1-a) y + \phi_{23}) \partial_y$ and the remaining term in $\partial_z$ is then the value of $d_3(\varphi \bullet X_{g_0} + \psi D_a \bullet P \omega)$. \\

Let $y_1 = \varphi \in H^0$, $y_2 = X_f \wedge \partial_z,\ y_3 = X_g \wedge \partial_z \in H^2$. The equation \eqref{formality-eq-k2} rewrites
\begin{equation*}
d_3(\varphi \bullet X_f \wedge \partial_z \bullet X_g \wedge \partial_z) = \delta_{CE}^0 \varphi_3(\varphi \bullet X_f \wedge \partial_z \bullet X_g \wedge \partial_z) + \varphi'_z \big(\{f_0,g_1\}+\{g_0,f_1\}\big) \partial_z.
\end{equation*}
The term in $\partial_z$ is not a coboundary, so we take it as the value of the $d_3$ term and set the value of the $\varphi_3$ term to zero. \\

Let $y_1 = \varphi \in H^0$, $y_2 = X_g \wedge \partial_z \in H^2$ and $y_3 = P \omega \in H^3$. The equation \eqref{formality-eq-k2} rewrites
\begin{equation*}
\begin{split}
d_3(& \varphi \bullet X_g \wedge \partial_z \bullet P\omega) \\
={}&  \delta_{CE}^1 \varphi_3(\varphi \bullet X_g \wedge \partial_z \bullet P\omega) + (P \varphi'_z X_{\partial_z g} - P \varphi''_z X_g - \phi_{13} \partial_x - \phi_{23} \partial_y) \wedge \partial_z \\
& - P \partial_z (g\varphi'_z)_1 \partial_x \wedge \partial_y + P \partial_y (g\varphi'_z)_1 \partial_x \wedge \partial_z - P \partial_x (g\varphi'_z)_1 \partial_y \wedge \partial_z,
\end{split}
\end{equation*}
with $\phi_{13},\phi_{23}$ coming from $\varphi_2(X_{g_0 \varphi'_z(0)},P\omega)$. Here the term in $\partial_x \wedge \partial_y$ is a coboundary and can be cancelled, the remaining expression is the value of $d_3$. \\

Let $y_1 = \varphi \in H^0$, $y_2 = P \omega,\ y_3 = Q \omega \in H^3$. The equation \eqref{formality-eq-k2} gives
\begin{equation*}
d_3(\varphi \bullet P \omega \bullet Q \omega) = \delta_{CE}^2 \varphi_3(\varphi \bullet P \omega \bullet Q \omega) - 2 PQ\varphi''_z \omega
\end{equation*}
We can cancel the coboundary $(- 2 PQ\varphi''_z)_1 \omega$ and it remains $d_3(\varphi \bullet P \omega \bullet Q \omega) = - 2 PQ\varphi''_z(0) \omega$. \\

Let $y_1 = X_{f_0} + \varphi D_a,\ y_2 = X_{g_0} + \psi D_a,\ y_3 = X_{h_0} + \xi D_a \in H^1$. The equation \eqref{formality-eq-k2} gives
\begin{equation*}
\begin{split}
d_3(& X_{f_0} + \varphi D_a,\bullet X_{g_0} + \psi D_a, \bullet X_{h_0} + \xi D_a) = - \xi'_z \big(D_a((f_0,\varphi,g_0,\psi)) - (f_0,\varphi,g_0,\psi)\big) \\
& + \psi'_z \big(D_a((f_0,\varphi,h_0,\xi)) - (f_0,\varphi,h_0,\xi)\big) - \varphi'_z \big(D_a((g_0,\psi,h_0,\xi)) - (g_0,\psi,h_0,\xi)\big) \\
& + (\varphi \psi'_z - \varphi'_z \psi + 2z (\varphi \psi''_z - \varphi''_z \psi))(D_a(h_0)-h_0) - (\varphi \xi'_z - \varphi'_z \xi - 2z (\varphi \xi''_z - \varphi''_z \xi))(D_a(g_0)-g_0) \\
& + (\psi \xi'_z - \psi'_z \xi + 2z (\psi \xi''_z - \psi''_z \xi))(D_a(f_0)-f_0). \\
\end{split}
\end{equation*}

The other cases also give very long expressions when plugged in equation \eqref{formality-eq-k2}.
\begin{gather*}
y_1 = X_{f_0} + \varphi D_a,\quad y_2 = X_{g_0} + \psi D_a \in H^1 \qquad y_3 = X_h \wedge \partial_z \in H^2 \\
y_1 = X_{f_0} + \varphi D_a,\quad y_2 = X_{g_0} + \psi D_a \in H^1 \qquad y_3 = P \omega \in H^3 \\
y_1 = X_{f_0} + \varphi D_a \in H^1,\qquad y_2 = X_g \wedge \partial_z,\quad y_3 = X_h \wedge \partial_z \in H^2 \\
y_1 = X_{f_0} + \varphi D_a \in H^1,\qquad y_2 = X_g \wedge \partial_z \in H^2,\qquad y_3 = P \omega \in H^3 \\
y_1 = X_f \wedge \partial_z,\quad y_2 = X_g \wedge \partial_z,\quad y_3 = X_h \wedge \partial_z \in H^2
\end{gather*}
Each time, the map $\varphi_3$ can be used to cancel as many term that is possible to obtain with some coboundary, but the remaining term for $d_3$ are still big expressions.

In the same way, we could obtain explicit values for higher maps $\varphi_k$ and $d_k$. Since there are several possibilities to choose the degrees of the elements to land in the non-trivial part of the complex $C[2]$, the $L_\infty$ structure obtained by solving inductively the formality equations has non-zero terms in each order.

\bibliographystyle{amsalpha}

\nocite{*}
\bibliography{biblio-heisenberg}

\newcommand{\etalchar}[1]{$^{#1}$}
\providecommand{\bysame}{\leavevmode\hbox to3em{\hrulefill}\thinspace}
\providecommand{\MR}{\relax\ifhmode\unskip\space\fi MR }
% \MRhref is called by the amsart/book/proc definition of \MR.
\providecommand{\MRhref}[2]{%
  \href{http://www.ams.org/mathscinet-getitem?mr=#1}{#2}
}
\providecommand{\href}[2]{#2}
\begin{thebibliography}{BGH{\etalchar{+}}05}

\bibitem[BFF{\etalchar{+}}78]{BFFLS78}
Fran\c{}cois Bayen, Mosh\'{e} Flato, Christian Fr{\o}nsdal, Andr\'{e}
  Lichnerowicz, and Daniel Sternheimer, \emph{Deformation theory and
  quantization {I}/ {II}}, Annals of Physics \textbf{111} (1978), no.~1,
  61--110 ; 111--151, \doi{10.1016/0003-4916(78)90224-5} and
  \doi{10.1016/0003-4916(78)90225-7}.

\bibitem[BGH{\etalchar{+}}05]{BGHHW05}
Martin Bordemann, Gr\'{e}gory Ginot, Gilles Halbout, Hans-Christian Herbig, and
  Stefan Waldmann, \emph{Formalit\'{e} ${G}_\infty$ adapt\'{e}e et
  star-repr\'{e}sentations sur des sous-vari\'{e}t\'{e}s co\"{i}sotropes},
  \href{http://arxiv.org/pdf/math/0504276.pdf}{arXiv : math/0504276v1
  [math.QA]}, 2005.

\bibitem[BM08]{BM08}
Martin Bordemann and Abdenacer Makhlouf,
  \emph{\href{http://www.springerlink.com/content/42w26ml594256t4r/fulltext.pdf}{Formality
  and Deformations of Universal Enveloping Algebras}}, International Journal of
  Theoretical Physics \textbf{47} (2008), 311--332.

\bibitem[BMP05]{BMP05}
Martin Bordemann, Abdenacer Makhlouf, and Toukaiddine Petit,
  \emph{\href{http://www.sciencedirect.com/science/article/pii/S0021869304005563}{D\'{e}formation
  par quantification et rigidit\'{e} des alg\`{e}bres enveloppantes}}, Journal
  of Algebra \textbf{285} (2005), no.~2, 623--648.

\bibitem[Bor]{B08}
Martin Bordemann, \emph{Deformation of a differential coderivation and
  ${L}_\infty$ structures}, Notes non publi\'{e}es.

\bibitem[Bro65]{Br65}
Ronald Brown, \emph{The twisted eilenberg-zilber theorem}, Edizioni Oderisi
  (Gubbio) (Simposio di~Topologia, ed.), Celebrazioni {A}rchimedee del {S}ecolo
  {XX}, (Messina, 1964), 1965, pp.~33--37.

\bibitem[CE56]{CE56}
Henri Cartan and Samuel Eilenberg,
  \emph{\href{http://www.archive.org/download/homologicalalgeb033541mbp/homologicalalgeb033541mbp.pdf}{Homological
  algebra}}, Princeton University Press, 1956.

\bibitem[Elc12]{El12}
Olivier Elchinger, \emph{\href{http://www.theses.fr/2012MULH4079}{Formality
  related to universal enveloping algebras and study of {H}om-(co){P}oisson
  algebras}}, Ph.{D}. thesis, Universit\'e de Haute-Alsace, Mulhouse, November
  2012.

\bibitem[Hue10]{Hueb10}
Johannes Huebschmann, \emph{The sh-{L}ie algebra perturbation {L}emma}, Forum
  Mathematicum \textbf{23} (2010), no.~4, 669--691,
  \href{http://arxiv.org/pdf/0710.2070v3.pdf}{arXiv:0710.2070 [math.AG]}.

\bibitem[Hue11]{Hueb11}
\bysame, \emph{The {L}ie {A}lgebra {P}erturbation {L}emma}, {H}igher
  {S}tructures in {G}eometry and {P}hysics (Alberto~S. Cattaneo, Anthony
  Giaquinto, and Ping Xu, eds.), Progress in Mathematics, vol. 287,
  Birkh\"{a}user Boston, 2011,
  \href{http://arxiv.org/pdf/0708.3977v2.pdf}{arXiv:0708.3977 [math.AG]},
  pp.~159--179 (English).

\bibitem[Kon03]{K03}
Maxim Kontsevitch,
  \emph{\href{http://www.springerlink.com/content/x0h825521k418w4x/fulltext.pdf}{Deformation
  quantization of {P}oisson manifolds, {I}}}, Letters of Mathematical Physics
  \textbf{66} (2003), no.~3, 157--216.

\end{thebibliography}

\end{document}